\documentclass[12pt]{article}
\font\smallit=cmti12
\usepackage{theorem,amssymb,amsmath}
\theorembodyfont{\rmfamily}

\newtheorem{thm}{Theorem}[section]
\newtheorem{lemma}[thm]{Lemma}
\newtheorem{cor}[thm]{Corollary}
\newtheorem{prop}[thm]{Proposition}
\newtheorem{defn}[thm]{Definition}

\begin{document}

\date{}
\title{\bf \large Avoiding Monochromatic Sequences With Special Gaps }
\author{Bruce M. Landman\\
{\smallit Department of Mathematics,
State University of West Georgia, Carrollton, GA 30118} 
\and Aaron Robertson \\
{\smallit Department of Mathematics,
Colgate University, Hamilton, NY 13346}
}
\maketitle

\begin{abstract}
\noindent
For $S$ a set of positive integers, and $k$ and $r$ fixed
positive integers, denote by $f(S,k;r)$ the least positive
integer $n$ (if it exists) such that within every $r$-coloring of
$\{1,2,...,n\}$ there must be a monochromatic sequence
$\{x_{1},x_{2},...,x_{k}\}$ with $x_{i}-x_{i-1} \in S$ for $2
\leq i \leq k$. We consider the existence of $f(S,k;r)$ for
various choices of $S$, as well as upper and lower bounds on this 
function. In particular, we show that this function
 exists for all $k$ if $S$ is an odd translate of the set of primes and $r=2$.

\end{abstract}

\baselineskip=14pt

\section{Introduction}

Van der Waerden's theorem on arithmetic progressions [9] states that for
every partition of the natural numbers {$\mathbb{N}$}  into $r$ sets, at
least one of the sets will contain arbitrarily long arithmetic progressions.
An equivalent form of this theorem says that for all positive integers
$k$ and $r$, there exists a positive integer $n=w(k;r)$ such that within
every $r$-coloring of $[1,n]=\{1,2,...,n\}$ there must be a monochromatic
$k$-term arithmetic progression. By replacing the set of arithmetic
progressions, $AP$, 
with another family ${\cal F}$ of sets,
one may ask if the corresponding theorem holds, i.e., is it true that
for all $k$ and $r$, there exists a positive integer $n=f(k,r)$ such that
for every $r$-coloring of $[1,n]$, there is a monochromatic $k$-term
member of ${\cal F}$? Examples may be found in [4,5,6,7,8].

In [5], the authors considered replacing $AP$ with a smaller collection of sets,
namely the collection of those arithmetic progressions $\{x+id: 0 \leq i \leq k-1\}$
whose common
 differences, $d$, belong to some prescribed set. Specifically, for a positive
  integer $r$, and $A$ a set of positive integers, call $A$ an {\em
$r$-large} set if for every $r$-coloring of the positive integers there
exist arbitrarily long monochromatic arithmetic progressions whose common
differences belong to $A$. Further, define a set to be {\em large}
if it is $r$-large for every $r$. They gave several sufficient conditions
and some necessary conditions for largeness and 2-largeness. They also
conjectured that any set that is 2-large must also be large.

In this paper we consider a property related to largeness. As with
largeness, we consider sequences where the differences between consecutive
terms belong to a prescribed set $S$; however, we do not insist that
the sequence be an arithmetic progression. We begin with the following
notation and definitions.

\noindent
{\bf Notation.}
For any string $w$ and any $t \in \mathbb{N}$, we denote by $w^t$ the
string $\underbrace{ww\cdots w}_{t}$.

\begin{defn}  Let $S \subseteq {\mathbb{N}}$. A sequence of positive
integers
 $\{x_{1},...,x_{k}\}$ is a $k$-term {\em $S$-diffsequence} if
 $x_{i}-x_{i-1} \in S$ for $2 \leq i \leq k$.
\end{defn}

\begin{defn} Let $r \in {\mathbb{N}}$. A set of positive integers $S$ is
called {\em $r$-accessible} is whenever {$\mathbb{N}$} is $r$-colored, there are
arbitrarily long monochromatic  $S$-diffsequences.
\end{defn}

\begin{defn} S is called {\em accessible} if $S$ is $r$-accessible
for all positive integers $r$.
\end{defn}

\begin{defn} If $S$ is not accessible, the {\em degree of
accessibility} of $S$ is the largest value of $r$ such that $S$ is
$r$-accessible. We denote this by DA($S$).
\end{defn}

 We denote by $f(S,k;r)$ the least positive integer $n$ (if it exists)
 such that for every $r$-coloring of $[1,n]$ there
 is a monochromatic $k$-term $S$-diffsequence.  Obviously, if $S \subseteq T$,
 then $f(S,k;r) \geq f(T,k;r)$.

 Denote the family of all accessible sets by ${\cal A}$ and the family
 of all $r$-accessible sets by ${\cal A}_{r}$. Likewise denote the families of
 large sets and $r$-large sets by ${\cal L}$ and ${\cal L}_{r}$, respectively.
 Clearly, ${\cal L} \subseteq {\cal A}$ and ${\cal L}_{r} \subseteq {\cal A}_{r}$
 for all $r$. In [5], it was conjectured that ${\cal L} = {\cal L}_{2}$.
 As we shall see, ${\cal A} \neq {\cal A}_{2}$ and
  ${\cal A}_{2} \neq {\cal L}_{2}$.  We still do not know whether
   ${\cal A} = {\cal L}$.

  In Section 2 we present some basic lemmas and
consider a few elementary examples.
Section 3 deals with sets consisting of certain congruence classes;
in particular, we will see that for each positive integer $d$, there is
some set having
$d$ as its degree of accessibility (this is in contrast to what
has been conjectured about large sets).  In Section 4 we prove
that for each odd positive integer $t$
   there are arbitrarily long sequences of primes
  $p_{1}<p_2< \cdots <p_{k}$ such that
   $p_{i}-p_{i-1} \in P+t$ for $2 \leq i \leq k$,
  where $P$ is the set of primes. From this it will follow that $P+t
\in \mathcal{A}_2$.
  Section 5 contains some open questions, as well as a table of
computer-generated
  values of $f(S,k;2)$ for several different sets $S$ and values $k$.

\section{A Few Simple Examples}

We begin with two  useful lemmas.

\begin{lemma}
Let $c \geq 0$ and $r \geq 2$, and let $S$ be a set of positive integers.
If every $(r-1)$-coloring of $S$ yields arbitrarily long monochromatic
$(S+c)$-diffsequences, then $S+c \in \mathcal{A}_r$.
\end{lemma}
{\em Proof.} Let $S = \{s_{i}: i \in {\mathbb{N}}\}$ and assume every
$(r-1)$-coloring of $S$ admits arbitrarily long monochromatic $(S+c)$-diffsequences.
Let $\chi$ be an $r$-coloring of ${\mathbb{N}}$. By induction on $k$, we
show that, under $\chi$, for all $k$ there are $k$-term monochromatic
$(S+c)$-diffsequences. Since there are obviously 1-term sequences, assume
$k \geq 1$ and that under $\chi$ there is a monochromatic
$(S+c)$-diffsequence
 $X=\{x_{1},...,x_{k}\}$. We may assume $X$ has the color red.
 Consider $A=\{x_{k}+ s_{i}+c: s_{i} \in S\}$. If some member of $A$ is
 colored red, then we have a red $(k+1)$-term $(S+c)$-diffsequence. Otherwise
 we have an $(r-1)$-coloring of $A$ and therefore, by the hypothesis, $A$ must
 contain arbitrarily long monochromatic $(S+c)$-diffsequences.
\hfill $\Box$

\noindent
{\em Remark 1}. The converse of Lemma 2.1 is false. As one example, let
$S = \{2\} \cup (2{\mathbb{N}} -1)$. Let $\chi$ be the 2-coloring
of
$S$ defined by $\chi(x) = 1$ if
 $x \equiv 1$(mod 4) or $x=2$, and  $\chi(x)=0$ if $x \equiv 3$(mod 4).
 Then $\chi$ does not yield arbitrarily long monochromatic
 $S$-diffsequences (there are none of length four).
 On the other hand, $S \in \mathcal{A}_3$ [8, Remark (5)], and in fact
 $f(S,k;3) \leq 6k^{2}-13k+6$; more generally, from
this same reference
it follows that if $m$ is even,
 and $j$ is a positive
 integer, then the set $\{jm\} \cup \{x: x \equiv \frac{m}{2}($mod $m)\}$
 is 3-accessible.

\begin{lemma}
Let $S$ be a set of positive integers and let $k,r,j \in {\mathbb{N}}$. If
 $f(S,k;r)=M$, then $f(jS,k;r) = j(M-1)+1$.
\end{lemma}
{\em Proof.}
Since $f(S,k;r)= M$, under any $r$-coloring of the set
$\{1,j+1,2j+1,..., (M-1)j+1\}$, there must exist a monochromatic $k$-term
$jS$-diffsequence.

On the other hand, let $\chi$ be an $r$-coloring of $[1,M-1]$ that avoids
monochromatic $S$-diffsequences of length $k$. Define
the $r$-coloring $\chi'$ of $[1,j(M-1)]$ by $\chi'[(i-1)j+1,ij] = \chi(i)$
for $i=1,2,...,M-1$. Assume, by way of contradiction, that $\chi'(x_{1}')=
\cdots = \chi'(x_{k}')$ with $x_{i}' - x_{i-1}'\in jS$ for $2 \leq i \leq k$.
Then, by the way $\chi'$ is defined, there exist $x_{1},...,x_{k}$,
monochromatic under $\chi$, belonging to $[1,M-1]$, with $x_{i}-x_{i-1} \in S$
for $2 \leq i \leq k$, a contradiction.
\hfill $\Box$

Using Lemma 2.1, with $r=2$, it is clear that the set $\{2^{i}: i \geq 0\}$
is 2-accessible. The following result tells us more.

\begin{thm}
Let $a \in {\mathbb{N}} \setminus \{1,3\}$, and let
\[S = \{(a-1)a^{j}: j=0,1,2,...\}  \cup \{(a-1)^{2}a^{j}:
j=0,1,2,...\}.\] Then $2 \leq DA(S) \leq a$. Furthermore, $f(S,k;2) \leq
a^{k}-a+1$ for all $k \geq 1$.
\end{thm}
{\em Proof.} To show that DA$(S) \leq a$, we exhibit an $(a+1)$-coloring
of  $\mathbb{N}$ which avoids monochromatic 2-term $S$-diffsequences.
Define
$\chi:  {\mathbb{N}} \rightarrow \{0,1,...,a\}$ by $\chi(x)=i$ where
$x \equiv i($mod $(a+1))$. Assume that $\chi(y)=\chi(z)$ and that $z-y \in S$.
By the definition of  $\chi$, $a+1$ divides $z-y$, and therefore either
$(a+1)|(a-1)a^{j}$ or $(a+1)|(a-1)^{2}a^{j}$ for some $j \geq 0$. Since
gcd$(a-1,a)=1$, we have
$(a+1)|a^{j}$ or $(a+1)|(a-1)^{2}$, but since $a \neq 3$, neither of these
is possible.

Now let $\alpha:[1,a^{k}-a+1] \rightarrow \{0,1\}$. To complete the proof
we show that under $\alpha$ there must be a monochromatic $k$-term
$S$-diffsequence. We do this by induction on $k$. Obviously, it holds
for $k=1$. Now assume $k \geq 2$, and that it holds for $k-1$. Let
$X =\{x_{1},...,x_{k-1}\}$ be a monochromatic $S$-diffsequence, say of
color $0$, that is contained in $[1,a^{k-1}-a+1]$.
 Consider the set $A=\{x_{k-1}+(a-1)a^{i}: i = 0,...,k-1\}$. Note that
 $A \subseteq [1,a^{k}-a+1]$. If
there exists $y \in A$ of color 0, then $X \cup \{y\}$ is a monochromatic
$k$-term $S$-diffsequence. If, on the other hand, no such $y$ exists, then
$A$ is a monochromatic $k$-term $S$-diffsequence.
\hfill $\Box$

\begin{cor}
If $S = \{2^{i}: i \geq 0\}$, then DA$(S)=2$ and
\[ 8(k-3)+1 \leq f(S,k,2) \leq 2^{k}-1     \]
for all $k \geq 3$.
\end{cor}
{\em Proof.} The fact that DA$(S)=2$ and the upper bound are immediate from
Theorem 2.3.

For the lower bound, first note that by direct calculation we find that
$f(S,3;2) = 7$ and $f(S,4;2)=11$. To complete the proof we show by
induction on $k$ that, for $k \geq 5$, the 2-coloring 
$\chi_{k} = (10010110)^{k-3}$ avoids
monochromatic $k$-term $S$-diffsequences. It is easy to check
directly that this statement is satisfied by
$k=5$. So now assume
$k \geq 5$, that
$\chi_{k}$ avoids $k$-term $S$-diffsequences, and consider $\chi_{k+1}$.

Let $X = \{x_{1},x_{2},...,x_{m}\}$ be a maximal length monochromatic
$S$-diffsequence under $\chi_{k+1}$.
We wish to show that $m \leq k$. Assume, by way of contradiction, that
$m \geq k+1$. Then $x_{m-1},x_{m} \in [8(k-3)+1,8(k-2)]$, or else the
inductive assumption would be contradicted. We consider the following cases.

\noindent
{\tt Case 1.} $\chi_{k+1}(X) =1$.

{\tt Subcase (a).} $x_{m-2} \in [8(k-3)+1,8(k-2)]$. 

In this subcase
we must have $x_{m-2}=8k-20$, $x_{m-1}=8k-18$, and $x_{m}=8k-17$. By the
structure of $\chi_{k}$, we see that $x_{m-3} \equiv 4$(mod 8). Hence,
there exists, under $\chi_{k}$, 
a monochromatic $S$-diffsequence of length $m-1$, contradicting
our assumption about $\chi_{k}$.

{\tt Subcase (b).} $x_{m-2} \not \in [8(k-3)+1,8(k-2)]$.

 Then $x_{m-1}
=8k-18$ and $x_{m}=8k-17$. By the structure of $\chi_{k}$, this implies
$x_{m-2} \equiv 6$(mod 8).  Then there is an $(m-1)$-term monochromatic 
$S$-diffsequence under $\chi_{k}$, a contradiction.

\noindent
{\tt Case 2.} $\chi_{k+1}(X)=0$.

{\tt Subcase (a).} $x_{m-2} \in [8(k-3)+1,8(k-2)]$.

 For this case we
have
$x_{m-2} = 8k-22$, $x_{m-1} = 8k-21$, and $x_{m} = 8k-19$. Then either
$x_{m-3} = 8(k-3)$ or $x_{m-3} \equiv 2$(mod 8).  If $x_{m-3} = 8(k-3)$,
then $m-3 \leq k-3$ because there can be only one term of an $S$-diffsequence
per 10010110-string, a contradiction. If $x_{m-3} \equiv 2$(mod 8), then there
is an $(m-1)$-term $S$-diffsequence of color 0 under $\chi_{k}$, a contradiction.

{\tt Subcase (b).} $x_{m-2} \not \in [8(k-3)+1,8(k-2)]$.

  Then
$x_{m-1}=8k-21$ and $x_{m}=8k-19$, and hence $x_{m-2} \equiv 3$(mod 8).
This is not possible, since there would then be a monochromatic
 $(m-1)$-term $S$-diffsequence under $\chi_{k}$.
\hfill $\Box$.

We next show that $a=2$ is the only value of $a$ for which $\{a^{i}:i\geq
0\}\in \mathcal{A}_2$.  To this end, we first prove the following lemma.

\begin{lemma}
 Let $m \geq 2$ and $i \geq 1$ with gcd$(i,m)=1$. Let
$S= \{x \in {\mathbb{N}}: x \equiv i$(mod $m)\}$. Then $S \not
\in \mathcal{A}_2$.
\end{lemma}
{\em Proof.}
Let $\chi:  {\mathbb{N}} \rightarrow \{0,1\}$ be defined by $\chi(x)=0$
if and only if $m$ divides $x$. Then since any $m$-term  $S$-diffsequence
must include some multiple of $m$ and some non-multiple of $m$, there is
no monochromatic $m$-term $S$-diffsequence.
\hfill $\Box$

\begin{prop}
If $a \geq 3$, then $\{a^{i}: i \geq 0\} \not \in \mathcal{A}_2$.
\end{prop}
{\em Proof.}  Let $T=\{a^{i}: i \geq 0\}$. Then
$T \subseteq \{x: x \equiv 1$(mod $a-1)\}$, and the result follows from
Lemma 2.5.
\hfill $\Box$

 In [5]
it was shown that if $A \not \in \mathcal{L}_r$ and $B
\not \in \mathcal{L}_s$, then
$A \cup B \not \in \mathcal{L}_{rs}$ (hence, whenever a finite union of
sets is large, at least one of the sets must be large). Essentially the
same proof can be used to prove the following lemma. We omit the proof.

\begin{lemma} If
$S \not \in \mathcal{A}_r$ and $T
\not \in \mathcal{A}_s$, then
$S \cup T \not \in \mathcal{A}_{rs}$.
\end{lemma}

It is easy to see, using Lemma 2.1, that the set $S=\{2\} \cup
(2{\mathbb{N}}-1)$ is 2-accessible, since the set of odd numbers itself is
an
$S$-diffsequence. The next theorem tells us more about $S$.

\begin{thm}
If $S=\{2\} \cup (2{\mathbb{N}}-1)$, then DA$(S)=3$. Furthermore,
$f(S,k;3) \leq 6k^{2}-13k+6$ and

\noindent
\begin{equation} f(S,k;2) = \left\{ \begin{array}{ll}
                        3k-4  & \mbox{if $k$ is odd} \\
                        3k-3  & \mbox{if $k$ is even}
                        \end{array}
                 \right. 
\end{equation}
\end{thm}
{\em Proof.}  The fact that $DA(S) \geq 3$ and
the bound for $f(S,k;3)$ were mentioned in Remark 1.
The fact that DA$(S) < 4$ follows from Lemma 2.7.  To see
this, note that the $2$-coloring of $\mathbb{N}$ given
by $001100110011 \ldots$ shows that $\{2\}$ is not
$2$-accessible and that the $2$-coloring of
$\mathbb{N}$ given by $01010101\ldots$ shows that
$2 \mathbb{N} - 1$ is not $2$-accessible.  Hence, we
have that DA$(S) = 3$.

 Let $g(k)$ be the function on the right side of (1).
We next show that $g(k)$ is an upper bound for $f(S,k;2)$.
By direct computation
it is easily checked that $f(S,k;2) \leq g(k)$ for $k=2$ and $k=3$. To
show this inequality holds for $k \geq 4$, it suffices to show
that for every $\{0,1\}$-coloring of $[1,g(k)]$, there exist $S$-diffsequences
$X_{1}=\{x_{1},x_{2},...,x_{k_{1}}\}$  and $X_{2}=\{y_{1},y_{2},...,y_{k_{2}}\}$
where $X_{1}$ has color 0, $X_{2}$ has color 1, and $k_{1}+k_{2} \geq 2k-1$.
This last fact is true for $k=4$ and $k=5$ by direct computation. To show
it holds for all $k$ we proceed by induction on $k$, showing that its truth
for $k$ implies its truth for $k+2$.

Assume that $k \geq 4$, and that for every 2-coloring of $[1,g(k)]$ there
exist monochromatic sequences $X_{1}$ and $X_{2}$ as described above.
Now 2-color $[1,g(k+2)] = [1,g(k)+6]$. To complete the proof we show
that there exists a $k'_{1}$-term $S$-diffsequence of color 0 and
a $k'_{2}$-term $S$-diffsequence of color 1 with
\begin{equation}
k'_{1}+k'_{2} \geq 2k+3.
\end{equation}
We assume, without loss of generality, that $k_{1} \geq k_{2}$.
Let $Y=\{x_{k_{1}}+1,x_{k_{1}}+2,...,x_{k_{1}}+6\}$.
We consider the following cases.

\noindent
{\tt Case 1.} There exist at least four elements of $Y$ that have
color 0. 

It is easy to see that these four  elements
may be appended to  $X_{1}$  to form a monochromatic
 $S$-diffsequence, and hence (2) holds.

\noindent
{\tt Case 2.} Exactly three elements of $Y$ have color 0.

Then
there exist two elements, $a$ and $b$, of these three  such that $X_{1} \cup
\{a,b\}$ forms a $(k_{1}+2)$-term $S$-diffsequence. Likewise there exist two
members, $c$ and $d$, of $Y$, having color 1 and such
that $X_{2} \cup \{c,d\}$ forms a $(k_{2}+2)$-term $S$-diffsequence.
This implies (2) for this case.

\noindent
{\tt Case 3.} Two or fewer elements of $Y$ have color 0.

 Then
we may extend $X_{2}$ to an $S$-diffsequence, monochromatic with color 1,
of length $k'_{2} \geq k_{2}+4$. Again (2) holds.

To complete the proof of the theorem, we show that
$f(S,k,2) \geq g(k)$ by exhibiting a 2-coloring of
 $[1,g(k)-1]$ that avoids monochromatic $k$-term $S$-diffsequences. We begin
 with the case in which $k$ is even. Let $C_{k}$ be the following
 coloring of $[1,3k-4]$:
 $ C_{k} =  1 (000111)^{\frac{k-2}{2}}0$. By symmetry
 it suffices to show there is no $k$-term $S$-diffsequence with color 1.
 We prove this by induction on $n$, where $k=2n$. Obviously the coloring
 10 avoids 2-term monochromatic $S$-diffsequences, and the coloring
 10001110 avoids 4-term monochromatic $S$-diffsequences, and hence the
 result holds for $n=1$ and $n=2$.

 Now assume $n \geq 2$, and that $C_{k}$ does not yield any $k$-term 
 monochromatic $S$-diffsequences with color 1. Now $C_{k+2} = C_{k}001110$.
 Let $X$ be a monochromatic $S$-diffsequence of color 1 in $C_{k}$ having
 maximal length. So $|X| < k$. Obviously, at least one of $\{3k-7,3k-6\}$ belongs to $X$.
 Hence $3k-5$ also belongs to $X$. Hence, at most two members of
 $\{3k-1,3k,3k+1\}$ may be tacked on to $X$ to form a monochromatic
 $S$-diffsequence. Thus, under the coloring $C_{k+2}$,
 there is no $k+2$-term $S$-diffsequence with color 1. This completes the
 proof for $k$ even.

Now consider the case in which $k$ is odd. Let $D_{k} = 11 (000111)^{\frac
{k-3}{2}}00$. The proof is
completed in a straightforward manner, similar to the even case,
 by induction on $n=(k-1)/2$, by showing
that the longest  $S$-diffsequence with color 1 cannot have length
greater than $k-1$. We omit the details.
\hfill $\Box$

It is a simple exercise to give an upper bound on $f(S,k;2)$ when $S$
is the set of Fibonacci numbers.  The proof is left to the reader.

 \begin{prop}
 Let $F=\{F_1,F_2,F_3,\dots\}=\{1,1,2,\dots\}$ be the sequence of
Fibonacci numbers. Then 
  $f(F,k,2) \leq F_{k+3}-2$.
 \end{prop}

 We conclude this section with the following simple result, which provides
 us with examples of very sparse sets which are nonetheless accessible.

 \begin{thm}
 Let $T \subseteq {\mathbb{N}}$ be infinite. Then 
 $T-T = \{t-s: s < t \mbox{ and $s,t \in T$} \} \in \mathcal{A}$.
\end{thm}
{\em Proof.}
Let $r \in {\mathbb{N}}$, and consider any $r$-coloring of $T-T$.  Fix
$s \in T$. Let $\{t_{1},t_{2},...\}
= \{t \in T: t > s\}$ where
$t_{1} < t_{2} <
\cdots$, and let $A=\{t_{i}-s: i =1,2,...\}$ . Obviously, there is some
color which contains an infinite subset, $B$, of $A$. 
Since $B$ is a $(T-T)$-diffsequence, by Lemma 2.1,
$T-T \in \mathcal{A}$.
\hfill $\Box$

 \section{Some Results on Sets of Congruence Classes}

 We now look at the accessibility of certain collections of congruence
 classes. In [5] it was proved that if a set $A$ belongs to
$\mathcal{L}_2$, then $A$ must contain a multiple of every positive
integer.
 We have seen  that this is not true if
 we replace $\mathcal{L}_2$ with $\mathcal{A}_2$ 
(see, for example, Corollary
2.4
 or Theorem 2.8). By the next lemma, we
 see that this condition is necessary in order for a set to be accessible.

 \begin{lemma}
 If $r \in {\mathbb{N}}$ and $S$ contains no multiple of $r$, then 
$S \not \in \mathcal{A}_r$.
 \end{lemma}
 {\em Proof.} Consider the $r$-coloring $\chi:{\mathbb{N}} \rightarrow
\{0,1,...,r-1\}$
 defined by $\chi(x) = i$ if $x \equiv i($mod $r$). This coloring avoids
 2-term monochromatic $S$-diffsequences.
 \hfill $\Box$

 We now consider
 the set of positive integers that, for a given $m$, are not multiples of $m$.
 We shall denote this set by $S_{m}$.  In [5], it was shown that $S_{m}
\not \in \mathcal{L}_2$,
and by Lemma 3.1, $S_{m} \not \in \mathcal{A}$. By the
following
 result, $S_{m} \in \mathcal{A}_2$ for $m > 2$, thus
giving another example
 for which ``2-accessible'' does not imply ``2-large.'' 

\begin{thm}
Let $m \geq 2$. Then DA$(S_{m})=m-1$.
\end{thm}
{\em Proof.} The fact that DA$(S_{m}) \leq m-1$ follows from Lemma 3.1.

To prove the reverse inequality, let $\chi$ be any $(m-2)$-coloring of $S$.
 Then some color
must contain an infinite number of elements from each of at least two of
the residue classes                      1 (mod $m)$, 2 (mod $m)$, 
..., $(m-1)$ (mod $m)$. Thus, some color contains arbitrarily
long $S$-diffsequences. By Lemma 2.1, $S
\in \mathcal{A}_{m-1}$, and the proof is complete.
\hfill $\Box$.

An immediate and noteworthy corollary of Theorem 3.2 is the following.

\begin{cor} Let $d\in \mathbb{N}$.  There exists 
$S \subseteq \mathbb{N}$ such that $DA(S)=d$.
\end{cor}

 In the
next theorem, we give the exact value of $f(S_{m},k;2)$ for $m=3$ and
$m=4$. We use $g(k)$ to
denote the right-hand side of (1) (from Theorem 2.8).

\begin{thm}
Let $k \geq 2$. Then
(i) $f(S_{3},k;2)=4k-5$  and (ii) $f(S_{4},k;2)=g(k)$.
\end{thm}
{\em Proof.} To prove $f(S_{3},k;2) \geq 4k-5$, consider the coloring
$\chi:[1,4k-6] \rightarrow \{0,1\}$, defined by $\chi(i) = 0$ if
$i \equiv 2$(mod 4) or $i \equiv 3$ (mod 4), and $\chi(i)=1$ if
$i \equiv 0$(mod 4) or $i \equiv 1$(mod 4).
In each color there are $k-2$ pairs of consecutive
elements that differ by 3. Hence in each color there are at most
$2(k-2)+1-(k-2)=k-1$ elements that can belong to the same $S_{3}$-diffsequence.
Hence $f(S_{3},k;2) > 4k-6$.

To prove the reverse inequality we will show the following stronger statement
is true: for every 2-coloring $\chi:[1,4k-5] \rightarrow \{0,1\}$
there exist $S_{3}$-diffsequences
$X=\{x_{1},x_{2},...,x_{k_{1}}\}$ and $Y=\{y_{1},y_{2},...,y_{k_{2}}\}$ with
$\chi(X)=0$ and $\chi(Y)=1$ and $k_{1}+k_{2} \geq 2k-1$.
We prove this last statement by induction on $k$. It is easy to check that
the statement holds for $k=2$.  Now assume $k \geq 2$, and that the result holds for
$k$. Let $\chi$ be any 2-coloring of $[1,4k-1]$. By inductive hypothesis, within
$[1,4k-5]$,
there exist monochromatic sequences $X$ and $Y$ as described above. Without
loss of generality, we assume $x_{k_{1}} \geq y_{k_{2}}$. We
consider three cases.

\noindent
{\tt Case 1.}  $x_{k_{1}} \equiv y_{k_{2}}$(mod 3). 

 Consider the
numbers $x_{k_{1}}+1$ and $x_{k_{1}}+2$. Regardless of their colors, we
now have monochromatic sets $\{x_{1},x_{2},...,x_{k'_{1}}\}$ and
$\{y_{1},y_{2},...,y_{k'_{2}}\}$ with $k'_{1}+k'_{2} = k_{1}+k_{2}+2 \geq 2k+1$.

\noindent
{\tt Case 2.} $x_{k_{1}} \equiv (y_{k_{2}}+1)$(mod 3).

 Let
$A = \{x_{k_{1}}+i: 1 \leq i \leq 4\}$, and let
 $A_{0}= \{x \in A: \chi(x)=0\}$ and $A_{1}=\{x \in A: \chi(x)=1\}$.
We may break this into  the following three subcases: (i) $A_{0}$ contains
 one of the pairs $\{x_{k_{1}}+1,x_{k_{1}}+2\}$,
 $\{x_{k_{1}}+1,x_{k_{1}}+3\}$, $\{x_{k_{1}}+2,x_{k_{1}}+3\}$,
$\{x_{k_{1}}+2,x_{k_{1}}+4\}$; (ii) $A_{1}$
 contains one of the pairs $\{x_{k_{1}}+1,x_{k_{1}}+2\}$, $\{x_{k_{1}}+1,
 x_{k_{1}}+3\}$, $\{x_{k_{1}}+3,x_{k_{1}}+4\}$; (iii)
 $A_{0}=\{x_{k_{1}}+1,x_{k_{2}}+4\}$ and $A_{1}=\{x_{k_{1}}+2,x_{k_{1}}+3\}$.
 In subcase (i), it is clear that there will be a $(k_{1}+2)$-term
 $S_{3}$-diffsequence with color 0, which gives the desired result. For
 subcase (ii), we have a $(k_{2}+2)$-term $S_{3}$-diffsequence with color 1.
 For subcase (iii) the monochromatic $S_{3}$-diffsequences
  $X \cup \{x_{k_{1}}+1\}$ and $Y \cup \{x_{k_{1}}+3\}$ yield the desired
  result.

  \noindent
{\tt Case 3.} $x_{k_{1}} \equiv (y_{k_{2}}+2)$(mod 3).

 Let $A$, $A_{0}$,
  and $A_{1}$ be defined as in Case 2. The following three subcases,
  which parallel the subcases of Case 2, yield the same respective results
  as those of Case 2:
\begin{itemize}
  \item[(i)] same as Case 2, subcase (i) 
   \item[(ii)] $A_{1}$ contains one of the pairs
$\{x_{k_{1}}+2,x_{k_{1}}+3\}$,
   $\{x_{k_{1}}+2,x_{k_{1}}+4\}$, $\{x_{k_{1}}+3,x_{k_{1}}+4\}$ 
   \item[(iii)] $A_{0} = \{x_{k_{1}}+3,x_{k_{1}+4}\}$ and
$A_{1}=\{x_{k_{1}}+1,
   x_{k_{1}}+2\}$.
\end{itemize}

   The fact that $f(S_{4},k;2) \leq g(k)$ follows immediately by
   Theorem 2.2, since $\{2\} \cup (2{\mathbb{N}}-1) \subseteq S_{4}$. 
   Also, the colorings $C_{k}$ and $D_{k}$ used in the proof of Theorem
2.8 not
   only avoid monochromatic
   $(\{2\} \cup (2{\mathbb{N}}-1))$-diffsequences, but they also
   avoid monochromatic $S_{4}$-diffsequences. Hence, $f(S_{4},k;2) \geq
g(k)$.
   \hfill $\Box$

   Although we do not have a formula for $f(S_{m},k;2)$ for $m > 4$,
   the next theorem gives a lower bound which we believe is the exact
   value of this function.

   \begin{thm}
   Let $m \geq 5$, and let $am \leq k < (a+1)m$.
    Then
   \[ 2k+2a-1 \leq f(S_{m},k;2).\]
   Furthermore, if $1 \leq k < m$, then $f(S_{m},k;2) = 2k-1$.
   \end{thm}
   {\em Proof.}
For the case in which $1 \leq k < m$,  choosing  any 2-coloring
of $[1,2k-2]$ such that there are $k-1$ elements of each color shows that
$f(S_{m},k,2) \geq 2k-1$. On the other hand, every 2-coloring of $[1,2k-1]$
yields a monochromatic set $\{x_{1},...,x_{k}\}$ such $x_{i}-x_{i-1} < m$
for $2 \leq i \leq k$.

Now let $a$ be as in the statement of the theorem. The lower bound follows
by observing that the 2-coloring
\[(10^{m-1})^{a}(1^{m-1}0)^{a}0^{k-a(m-1)-1}1^{k-a(m-1)-1} \]
avoids monochromatic $k$-term $S_{m}$-diffsequences.
\hfill $\Box$

  By Lemmas 2.2 and 2.5, we know that if 
$m \geq 2$ and $c \in \{1,2,...,m-1\}$
   then the set of positive
  integers that are congruent to $c$ (mod $m$) is not 2-accessible.
  We would like to know about the function $f(S,k;2)$ when $S$ is the
  union of more than one congruence class modulo $m$,
  other than the case in which $S = S_{m}$. We present one example, but
  thus far have not found a general result.

  \begin{prop}
  Let $S = \{x: 3 \nmid x$ and $4 \nmid  x\}$. Then
$f(S,k;2) = 7k-12$
  for $k \geq 3$.
  \end{prop}
  {\em Proof.} To see that $f(S,k;2) \geq 7k-12$ for $k$ even, note
  that the coloring $1(1001100 0110011)^{\frac{k-2}{2}}$ avoids monochromatic
  $k$-term $S$-diffsequences. For $k$ odd, the same property holds for
  the coloring $1(1001100 0110011)^{\frac{k-3}{2}}(1001100)$.

  To establish the upper bound we prove the following stronger statement:
  every 2-coloring $\chi$ of $[1,7k-12]$ admits monochromatic $S$-diffsequences
  $X = \{x_{1},x_{2},...,x_{m}\}$ and $Y=\{y_{1},y_{2},...y_{n}\}$, with
  $\chi(X)=0$ and $\chi(Y)=1$, such that $m+n \geq 2k-1$.
We prove this by induction on $k$. Direct calculation shows that
$f(S,3;2)=9$.  Now assume the statement holds
for $k$, and let $\chi$ be a 2-coloring of $[1,7k-5]$. Without loss of
generality, assume $x_{m} > y_{n}$. We may consider twelve cases, one
each for the congruence class,  modulo 12, that
contains $x_{m}-y_{n}$. We give the
details for two of these cases; the others are straightforward, and we omit
them.

\noindent
{\tt Case 5.} $x_{m}-y_{n} \equiv$ 5(mod 12).  

If $\chi(x_{m}+1) = 0$,
then regardless of the value of $\chi(x_{m}+2)$ we are done. Hence, assume
that $\chi(x_{m}+1)=1$.  If $\chi(x_{m}+4)=0$, then regardless of
the value of $\chi(x_{m}+2)$ we are done. Hence, assume $\chi(x_{m}+4)=0$.
Hence we may assume that $\chi(x_{m}+2)=0$, which implies that $\chi(x_{m}+3)=0$.
Now, regardless of the value of $\chi(x_{m}+7)$, we are done.

\noindent
{\tt Case 6.} $x_{m}-y_{n} \equiv$ 6(mod 12). 

 If $\chi(x_{m}+2)=0$, then
regardless of the value of $\chi(x_{m}+1)$, we are done. Hence, assume
that $\chi(x_{m}+2)=1$. This implies that we may assume that $\chi(x_{m}+1)=0$,
which in turn allows us to assume that $\chi(x_{m}+6)=1$. This implies
that $\chi(x_{m}+3)=1$, or else we are done. From this we may assume
that $\chi(x_{m}+4)=\chi(x_{m}+5)=0$. Now, if $\chi(x_{m}+7)=1$, we are
done, so assume that $\chi(x_{m}+7)=0$. Then we have that $X \cup
\{x_{m}+5,x_{m}+7\}$ is monochromatic, completing this case.
\hfill $\Box$

 \section{Translations Of The Set Of Primes}

 In [5], the question was raised as to whether there exist any translations
 of $P$, the set of primes, that are large, or for that matter 2-large.
 Since a 2-large set must contain a multiple of every integer, $P+e
\not \in \mathcal{L}_2$ if $e$ is even. Likewise, by Lemma 3.1, if $e$ is
even, then
 $P+e \not \in \mathcal{A}$, and $P$ itself is not
$4$-accessible.  In fact, $P \not \in \mathcal{A}_3$.
To see this, color the multiples of $9$ green, the
remaining even numbers red, and the remaining odd numbers blue.  It is
easy to see that any sequences of $9$ reds, $9$ blues, or $2$ greens
must have numbers which differ by a non-prime.  We do not
 know whether $P$ is 2-accessible or
 whether any even translation of $P$ is 2-accessible. On the other hand,
   as we shall see in this section, all
 odd translations of $P$ are 2-accessible.

We use an application, given as
Theorem 4.1 below, of a theorem due
to Balog [1].  Before stating the theorem,
we introduce some notation.

Let ${\bf
b}= (b_{1},b_{2},...,b_{k}) \in \mathbb{Z}^k$,
$p \in P$, 
and $x \in 
\mathbb{R}^+$. We define:

$\pi(x; {\bf b}) = |\{n:1<n+b_i \leq x \,\, \mathrm{is \,\, prime \,\,
for \,\, every \,\,}
 1 \leq i \leq k\}|$; 
 
$\rho(p) =\rho(p;{\bf b}) = |\{n \,(\mbox{mod $p$}):
(n+b_{1})(n+b_{2})...(n+b_{k}) \equiv 0 \,(\mbox{mod $p$})\}|;$

$ \sigma({\bf b}) = 
\left\{
\begin{array}{ll}
\prod_{p \in P} \left(1 - \frac{1}{p}\right)^{-k} \left(1 -
\frac{\rho(p)}{p}\right) &
\mathrm{if \,\,} \rho(p) < p \mathrm{\,\,for \,\, all \,\,
 primes \,\,} p \\ \\
0 & \mathrm{otherwise;}
\end{array}
\right.
$

$T(x;{\bf b}) =  \underset{1\leq i \leq k}
{\underset{1<n+b_i \leq x} \sum} \frac{1}{\log (n+b_{1}) \log
(n+b_{2})...\log (n+b_{k})}.$

Before stating Theorem 4.1, we remind the reader of the 
following notation.

\noindent
{\bf Notation.}  Let $f(x)$ and $g(x)$ be functions
and let $k$ be a parameter.
We write $f(x) \gg g(x)$ if there exists
a constant, $c$, such that
$\lim_{x \rightarrow \infty} \frac{f(x)}{g(x)} \geq c$.
We write $f(x) \gg_k g(x)$ if the
constant $c$ is dependent upon $k$.

\begin{thm} (Balog).  Let $k \in \mathbb{Z}^+$,
let $x \in {\mathbb R}^+$ be sufficiently large, let
$t$ be a fixed nonnegative integer, and let 
$$B = \{ 
(0,q_1+t,\dots,\sum_{i=1}^{k-1} (q_i+t)):
 q_i \in P,
k \leq q_i \leq x/2k, 1 \leq i\leq k-1 \}.
$$ Define
$Z = \{{\bf b}=(b_1,\dots,b_k) \in B: \{n:1<n+b_k \leq x\} \neq
\emptyset\}$.  Then
$$
\sum_{{\bf b} \in Z} |\pi(x;{\bf b}) - \sigma({\bf b})
T(x;{\bf b}) | \ll_k \frac{x^k}{\log^{2k} x}.
$$
\end{thm}

\noindent
{\em Remark.} This follows from Balog's theorem ([1], p.49) with $A=2k$,
$c=0$,
$D=1$, and
$a_i=1$ for
$i=1,2,\dots,k$, since
$B$ is a subset of $Z$ as defined in
Balog's theorem. 
\hfill{$\Box$}

We will need the following technical lemma.  Before stating the
lemma we give a definition.

\noindent
{\bf Definition 5.} Let $p$ be prime.  We call a set of polynomials 
${\mathcal P} \subseteq \mathbb{Z}[y]$ $p$-{\it admissible} if there
exists an integer $h$ such that
$p$ does not divide any element of the set
when $y=h$.  If  $\mathcal P$ is
$p$-admissible for all primes $p$, we
call $\mathcal P$ {\it admissible}.

\begin{lemma}
Let $k \geq 2$ and let $t \geq 1$ be odd.
For $(z_1,z_2, \dots, z_{k-1}) \in \mathbb{Z}^{k-1}$, define the set of
polynomials
$$
Y_{(z_1,z_2, \dots, z_{k-1})}(y) = \left\{
y+\sum_{j=1}^{i-1} (z_j+t): 1 \leq i \leq k\right\}
\subseteq \mathbb{Z}[y]
$$ 
and let
$$
M=\left\{(q_1,\dots,q_{k-1}): k<q_1,\dots,q_{k-1} \leq x/2k \mbox{ are 
primes and } Y_{(q_1,\dots,q_{k-1})}(y) \mbox{ is admissible}
\right\}
$$  
for $x \in \mathbb{R}^+$ sufficiently large.
Then $|M| \gg_k (\frac{x}{\log x})^{k-1}$.
\end{lemma}

\noindent
{\em Proof.} 
Our approach is to show that for ``most" 
$(k-1)$-tuples of primes, $Y_{(q_1,\dots,q_{k-1})}(y)$ is admissible.
First of all, for any $(k-1)$-tuple of primes
$(q_1,\dots,q_{k-1})$, it is clear that
$Y_{(q_1,\dots,q_{k-1})}(y)$
is $p$-admissible for any prime $p \geq k$. 
Hence, we need to consider the $q$-admissibility
for primes $q<k$.
To this end,
consider those primes $r_{1}=2$, $r_{2}=3,
\dots,r_{d}$ less than $k$.  We will obtain
a lower bound for the number of  $(k-1)$-tuples
of primes
which are $r_i$-admissible for all $1 \leq i \leq d$.

Let $h$ be odd.  Below, we will find $q_1,\dots,q_{k-1}$
such that $Y_{(q_1,\dots,q_{k-1})}(y)$ is admissible
with $y=h$.
(We are in fact proving something stronger:
we prove that $Y_{(q_1,\dots,q_{k-1})}(y)$ is $r_i$-admissible 
with $y=h$ for $1 \leq i \leq d$, i.e., the same $h$
works for all $r_i$.)
So that
 $r_{i} \nmid (h+q_{1}+t)$ for $1 \leq i
 \leq d$, it is sufficient that for each $i$,
 $q_1 \not \equiv -h-t \,\, (\mathrm{mod\,\,}r_i)$ 
Letting $m = \prod_{i=1}^{d} r_{i}$ we
need only have
$q_1$ belong to one specific residue class $c_1$ (mod $m$), where
$\gcd(c_1,m)=1$.  By Dirichlet's theorem for primes in arithmetic
progressions, we have
$\gg_k \frac{x}{\log x}$ choices for
$q_{1}$. 

Similarly, once $h,q_{1},q_{2},...,q_{j-1}$
have been chosen, we may choose $q_j$ so that for each $r_{i}$,
$q_{j}$ avoids one specific residue class modulo $r_{i}$. Hence we need
only choose $q_{j}$ so that it does not belong to any of the residue
classes
$-(h+q_{1}+q_{2}+ \cdots q_{j-1} + jt)$ (mod $r_{i}$), 
$1 \leq i \leq d$. So it suffices to
have
$q_{j}$ belong to one specific
congruence class $c_j$(mod $m$), with $\gcd(c_j,m)=1$. 

  Combining this criteria for all
primes less than $k$, we have at least $\prod_{i=2}^{d} (r_i-2)$
reduced residue classes modulo $\prod_{i=2}^{d} r_i$.
By Dirichlet's Theorem we have
$\gg_k \frac{x}{\log x}$ choices for each $q_i$, and thus
$\gg_k (\frac{x}{\log x})^{k-1}$ choices for the $(k-1)$-tuple
of primes $(q_1,q_2,\dots,q_{k-1})$ that belong
to $M$.
\hfill{$\Box$}

Using Theorem 4.1 and Lemma 4.2, we have
the following result.

\begin{lemma}
For $k \geq 2$, $t \geq 1$ and odd, and 
$x \in {\bf R}^+$, sufficiently large,
define 
$$
W =
\left\{(p,q_1,\dots,q_{k-1}): p,q_1,\dots,q_{k-1}
\mbox{ are primes and }
k<q_1,q_2,\dots,q_{k-1} \leq \frac{x}{2k} \right\}.
$$
For $1 \leq i \leq k-1$, let 
$$
S_i=\left\{(p,q_1,\dots,q_{k-1}) \in W
 : p+\sum_{j=1}^i (q_j+t) \leq x 
\mbox{ is
prime} \right\}
$$  and let $S =\bigcap_{i=1}^{k-1} S_i$.
Then $|S| \gg_{k} \frac{x^k}{\log^{2k-1}x}.$
\end{lemma}

\noindent
{\it Proof}.  
We  use the notation from Theorem 4.1 and Lemma 4.2; in particular,
${\bf b} =(0,q_1+t,\dots,\sum_{i=1}^{k-1} (q_i+t))$
and $M$ is as in Lemma 4.2.
In order to apply Theorem 4.1,
we first  
obtain effective bounds for
$\rho, \sigma$, and $T$. 

From Theorem 4.1 we see that we may restrict
our attention to those ${\bf b}$ such that
$\sigma({\bf b})>0$ and use the
same given bound (since this restriction
reduces the size of the sum in Theorem 4.1).
  
It is well known that $\sigma({\bf b}) < \infty$
(see [2], for example).  We next show
that for all  ${\bf b}=(q_1,\dots,q_{k-1}) \in M$
we have
$\sigma({\bf b})>0$.    
Since

\begin{enumerate}
\item For any $(q_1,\dots,q_{k-1}) \in M$
we have that
$Y(y) = \{ y + \sum_{j=1}^i (q_i+t): 1 \leq i \leq k-1\} \subseteq
\mathbb{Z}[y]$ is admissible,

\item[] and
  
\item
$Y(y)$  
is admissible
if and only if $\rho(p;{\bf b}) \leq p-1$
for each prime $p$, 
\end{enumerate} 
we see that for all ${\bf b}=(q_1,\dots,q_{k-1}) \in M$
we have $\rho(p;{\bf b}) < p$ for all primes $p$.
Since it is also true that $\rho(p;{\bf b}) \leq k$ for any
prime we have
\begin{equation}
\sigma({\bf b}) \geq  \prod_{p\leq k}
\left(1 - \frac{1}{p}\right)^{-k} \left(1 - \frac{p-1}{p}\right)
\prod_{p>k} \left(1 - \frac{1}{p}\right)^{-k} \left(1 - \frac{k}{p}\right)
=
\sigma_k,
\end{equation}
a constant dependent upon only $k$.    We now
show that $\sigma_k>0$.  

Clearly, we have
the finite product in (3) positive, so we must show that
the infinite product in (3) converges to a positive constant.
To this end, let $1+a_p=\left(1 - 1/p\right)^{-k} \left(1 -
k/p\right)$.  By the binomial theorem, we have 
$a_p =\frac{ -\sum_{i=2}^k (-1)^{k-i}{k \choose i} p^{-i}}
{(1 - 1/p)^k}$.   Since
$|a_p| \leq \frac{\sum_{i=2}^k {k \choose i} p^{-i}}
{(1 -1/p)^k} \leq \frac{\sum_{i=2}^k {k \choose i} p^{-2}}
{(1 - 1/p)^k} \leq \frac{\sum_{i=2}^k {k \choose i} p^{-2}}
{1/2^k} = 2^k(2^k-k-1)p^{-2}$, we see that
$\sum_{p \in P} a_p$ converges absolutely.
It follows that $\prod_{p \in P} (1+a_p)$ converges to a positive
number. 
Thus, from (3),
\begin{equation}
\mbox{for all } {\bf b} \in M,\,\,
\sigma({\bf b}) \geq 
\sigma_k > 0.
\end{equation}

We next bound $T(x;{\bf b})$ 
by using
\begin{equation*}
\begin{split}
|\{n:1<n+b_i \leq x, 1 \leq i \leq k \}|
&=(x-b_{k})+O(1)\\
&= x - \sum_{i=1}^{k-1} (q_i+t) + O(1) \\
&> x - \sum_{i=1}^{k-1} q_i -kt + O(1) \\
&> x - k\left(\frac{x}{2k} \right) + O(1) \\
&= \frac{x}{2} + O(1).
\end{split}
\end{equation*}
\vskip -40pt \hfill (5)
\newpage
\addtocounter{equation}{1}
\noindent
This gives us 
\begin{equation}
T(x;{\bf b}) >\left(\frac{x}{2}+O(1) \right) \frac{1}{\log^k
x}.
\end{equation}

From (5) we may
apply Theorem 4.1 to get
\begin{equation}
 \sum_{(q_1,\dots,q_{k-1}) \in M}
\Bigl| |\{n:n+b_i \,\, {\mathrm{
prime}}, 1 \leq i \leq k \}|  - \sigma({\bf b})T(x;{\bf b})  \Bigr| \ll_k
\frac{x^{k}}{\log^{2k}x} . 
\end{equation}

Using the bounds from (3), (4), and (6) along with
Lemma 4.2, inequality (7) yields

$$
\begin{array}{rl}
|S| &\geq  \underset{(q_1,\dots,q_{k-1}) \in M}{\sum} 
|\{n:n+b_i \,\, {\mathrm{
prime}}, 1 \leq i \leq k \}|\\
&\gg_k \underset{(q_1,\dots,q_{k-1}) \in M}{\sum} \sigma({\bf
b})T(x;{\bf b}) -O\left(\frac{x^k}{\log^{2k}x}\right)\\
&\gg_{k} \sigma_k |M|
\left(\frac{x}{2}+O(1) \right)\left(\frac{x}{2\log^k x}\right) - O
\left(\frac{x^k}{\log^{2k}x}\right)\\
&\gg_{k} \sigma_k \left(\frac{x}{\log x} \right)^{k-1}
\left(\frac{x}{2\log^k x}\right) - O
\left(\frac{x^k}{\log^{2k}x}\right)\\ &\gg_k \frac{x^k}{\log^{2k-1} x}\\
\end{array}
$$
for $x$ sufficiently large.
\hfill{$\Box$}

Using Lemma 4.3 we have the following result
concerning the existence of arbitrarily long
sequences of primes with ``special gaps."

\begin{thm}
Let $t \in \mathbb{N}$ be odd.   For any $k \geq 2$, there
exist
$p_{1},p_{2},...,p_{k} \in P$ such that $p_{i}-p_{i-1} \in P+t$ for
$i=2,\dots,k$.
\end{thm}

\noindent
{\em Proof.}  By Lemma 4.3 we may choose primes $p_1$, $q_1,
\dots,q_{k-1}$ so that 
$$p_i=p_1+\sum_{j=1}^{i-1} (q_j+t) \in P,$$
$2 \leq i \leq k-1$. Since $p_{i}-p_{i-1} = q_{i-1}+t$
for $i=2,\dots,k$, we are done.
\hfill{$\Box$}

\vskip 10pt
Combining Theorem 4.4 with Lemma 2.1, we have the following 
immediate corollary.

\begin{cor}
If $t$ is odd, then $P+t \in \mathcal{A}_2$.
\end{cor}

\section {Open Questions And Some Exact Values}

There are many interesting questions left unanswered about accessibility.
Here is a list of some that we would very much like to answer.

\noindent
1. True or false: ${\cal A}={\cal L}$? (this conjecture was posed by Tom Brown
[3]). It was proved in [5] that a set $S=\{s_{1},s_{2},...\}$ cannot be
large if $\liminf \frac{s_{i+1}}{s_{i}} > 1$. From the present paper we know
that $T-T \in \mathcal{A}$  for any infinite $T$, so that a
set can be very ``sparse'' and still be accessible. Perhaps an example
showing the answer to the above
question is false can be found by choosing the correct $T$; for example,
is the set $T-T$ large if $T = \{n!: n \in {\mathbb{N}}\}$?

\noindent
2. For $S = \{2^{i}: i \geq 0\}$, what is the exact value of $f(S,k;2)$?
We believe the lower bound of Corollary 2.1 is the exact value for $k
\geq 5$. In Table 1 (below) we give the first few values of this function.

\noindent
3. What is the exact value of $f(S,k;3)$ where 
$S = \{2\} \cup (2{\mathbb{N}}-1)$?

\noindent
4. What is a formula for $f(S_{m},k;2)$ that generalizes Theorem 3.4?
   Calculations for the case $m=6$ support the conjecture that the
   lower bound of Theorem 3.5 is the actual value of $f$, i.e., that for
   $k \geq 2$,
   \[ f(S_{6},k;2) = \left\{ \begin{array}{ll}
                      (5k-4)/2    &     \mbox{ if $k \equiv 2$(mod 4)} \\
                      (5k-5)/2 & \mbox{ if $k \equiv 3$(mod 4)} \\
                      (5k-6)/2 & \mbox{ if $k \equiv 0$(mod 4)}\\
                      (5k-7)/2 & \mbox{ if $k \equiv 1$(mod 4)}
                      \end{array}
                      \right.  \]

\noindent
5. If $t$ is an odd positive integer, what is DA($P+t)$? Moreover, is
it true that for every 2-coloring of $P$, there exist arbitrarily long
monochromatic $(P+t)$-diffsequences? If the answer to the latter question
is true, then by Lemma 2.1, $P+t \in \mathcal{A}_3$.

\noindent
6. What is the order of magnitude of $f(P+t,k;2)$ for a fixed odd positive
integer $t$? Table 1 below includes some specific values of this function.

\noindent
7. As stated earlier, $P \not \in \mathcal{A}_3$.
Is $P\in \mathcal{A}_2$? If so, what is the magnitude of
$f(P,k;2)$? We have calculated the first several values of $f(P,k;2)$
(see Table 1). 

\noindent
8. What is the degree of accessibility of the set of Fibonacci numbers?
   What is the order of magnitude of $f(F,k;2)$?

\noindent
9. What can we say about DA($S$) and $f(S,k;2)$ where $S$ is the union
of more than one congruence class modulo $m$? That is, generalize
Proposition 3.6.

The following table gives the exact value of $f(S,k;2)$ for various
choices of $S$ and $k$. The symbols $T$, $F$, and $P$ denote $\{2^{i}: i \geq 0\}$,
the set of Fibonacci
numbers, and the set of primes, respectively.

\begin{center}
\begin{tabular}{||r|r|r|r|r|r|r|r||} \hline\hline
\multicolumn{1}{||c|}{$S \setminus k$}      &       \multicolumn{1}{c|}{2}
   &  \multicolumn{1}{c|}{3}  &  \multicolumn{1}{c|}{4}  &
    \multicolumn{1}{c|}{5}  &  \multicolumn{1}{c|}{6}
    & \multicolumn{1}{c|}{7}  & \multicolumn{1}{c||}{8}\\ \hline
    $T$                    &  3  &  7  &  11  &  17  &  25  &  35&51\\ 
    $F$                    &  3  &  5  &   9  &  11  &  15  & 19  &  21 \\
    $P$                    &  5  &  9  &  13  &  21  &  25  & 33  &  ? \\
    $P+1$                  &  7  & 13  &  21  &  27  &  35  &  ?  &  ? \\
    $P+2$                  &  9  & 17  &  25  &  33  &   ?  &  ?  &  ? \\
    $P+3$                  & 11  & 21  &  31  &  42  &   ?  &  ?  &  ? \\
    $P+4$                  & 13  & 25  &  37  &   ?  &   ?  &  ?  &  ? \\
    $P+5$                  & 15  & 29  &   ?  &   ?  &   ?  &  ?  &  ? \\
    $P+6$                  & 17  & 33  &   ?  &   ?  &   ?  &  ?  &  ? \\
    $P+7$                  & 19  & 37  &   ?  &   ?  &   ?  &  ?  &  ? \\
    $S_{5}$                &  3  &  5  &   7  &  11  &  13  & 15  & 19 \\
    $S_{6}$                &  3  &  5  &   7  &   9  &  13  & 15  & 17\\
\hline\hline
    \end{tabular}
    \end{center}

\noindent
{\bf Acknowledgments.}  We would like to thank Andrew
Granville for guiding us to Balog's theorem and for
his invaluable assistance with the proof of Lemma 4.3.
We would also like to thank Scott Ahlgren for helping
with some details of the proof of Lemma 4.3.

 \end{document}